%% file: main_wordrepgraphs.tex
\documentclass[
a4paper,
envcountsame,
]{llncs}

\usepackage[utf8]{inputenc}
\usepackage[textsize=scriptsize,obeyFinal]{todonotes}

\usepackage[T1]{fontenc}
\usepackage{color}
\usepackage{xcolor}
\usepackage{graphicx}
\usepackage{amsmath}
\usepackage[parfill]{parskip}
\usepackage{xfrac}
\usepackage[]{easytable}
\usepackage{tikz}
\usepackage{mathtools}
\usepackage{alltt}
\usepackage{hyperref}
\usepackage[capitalise]{cleveref}
\usepackage{mathrsfs}
\usepackage{todonotes}
\usepackage{url}

\usepackage{amsmath, amsfonts, amssymb, dsfont}

\usepackage{enumerate}

\usepackage{wrapfig}
\usepackage{float}

\usetikzlibrary{positioning, shapes,patterns,shadows.blur, arrows,decorations,arrows,automata,shadows,patterns,chains,graphs,calc,intersections,matrix,fit,shapes,chains,decorations.pathreplacing}
\usepackage{paralist}
\usepackage{algorithm2e}

\usepackage{booktabs}
\usepackage{float}

\usepackage{caption}
\usepackage{subcaption}

\DeclareMathOperator{\alphabet}{alph}

\DeclareMathOperator{\loc}{loc}
\DeclareMathOperator{\Free}{Free}
\DeclareMathOperator{\cwd}{cwd}
\DeclareMathOperator{\lab}{label}
\DeclareMathOperator{\G}{G}

\newcommand{\m}{\setminus}

\newcommand{\ceil}[1]{\lceil #1 \rceil}

\newcommand{\N}{\mathbb{N}}
\newcommand{\R}{\mathbb{R}}

\newcommand{\E}{\mathscr{E}}
\renewcommand{\R}{\mathcal{R}}
\renewcommand{\L}{\mathcal{L}}
\renewcommand{\t}[1]{\mathtt{#1}}

\newcommand{\ta}{\mathtt a}
\newcommand{\tb}{\mathtt b}
\newcommand{\tc}{\mathtt c}
\newcommand{\td}{\mathtt d}
\newcommand{\te}{\mathtt e}

\newcommand{\tg}{\mathtt g}
\newcommand{\tl}{\mathtt l}
\newcommand{\tn}{\mathtt n}
\renewcommand{\to}{\mathtt o}
\newcommand{\tr}{\mathtt r}
\newcommand{\tp}{\mathtt p}

\renewcommand{\alph}{\alphabet}

\def\nth#1{#1$^{\text{th}}$}

\definecolor{cau}{RGB}{156,10,125} 

\newif\ifpaper
\papertrue 
\paperfalse 

\title{Word-Representable Graphs and Locality of Words}

\titlerunning{Word-Representable Graphs and Locality of Words}

\author{Philipp Böll$^1$\and Pamela Fleischmann$^2$\and Annika Huch$^3$\ Jana Kreiß$^4$\and Tim Löck$^5$\and\\Kajus Park$^6$\and Max Wiedenhöft$^7$}

\institute{Kiel University, Germany, \email{$\{$stu238462$^1$,stu201977$^4$,stu229765$^5$,stu222002$^6$$\}$@mail.uni-kiel.de, $\{$fpa$^2$,ahu$^3$,maw$^7$$\}$@informatik.uni-kiel.de}}

\authorrunning{P. Böll\and P. Fleischmann\and A. Huch\and J. Kreiß\and T. Löck\and K. Park\and M. Wiedenhöft}

\begin{document}

\maketitle

\begin{abstract}
In this work, we investigate the relationship between $k$-repre\-sentable graphs and graphs representable by $k$-local words.
In particular, we show that every graph representable by a $k$-local word is $(k+1)$-representable.
A previous result about graphs represented by $1$-local words is revisited with new insights.
Moreover, we investigate both classes of graphs w.r.t. their being hereditary and in particular the speed as a measure. We prove that the latter ones belong to the factorial layer and that the graphs in this classes have bounded clique-width.
\end{abstract}

\section{Introduction}
\input{intro}

\section{Preliminaries}\label{prelims}
\input{prelims}

\section{About $k$-Representability and Locality}\label{k+1}
\input{k+1.tex}

\section{Forbidden Subgraphs and Speed}\label{speed}
\input{speed.tex}

\section{Clique-Width}\label{cwd}
\input{cwd.tex}

\section{Conclusion}\label{conclusion}
\input{conclusion}

\bibliographystyle{plainurl}
\bibliography{refs}

%
  

\end{document}

%% file: intro.tex
The locality of a word is a complexity measure
introduced by Day et al. in \cite{LocalPatterns}
for pattern matching.
A pattern is a word consisting of letters (constants) and variables.
It is matched to a word by replacing the variables with letters.
The matching problem,
which is about deciding if a pattern can be matched to a certain word,
is NP-complete in general (cf. \cite{DBLP:journals/jcss/Angluin80}),
but can be solved in polynomial time for patterns with bounded locality.
The idea of locality is to iteratively mark all letters in a word according to a given marking sequence, which is an
enumeration of the word's letters.
In each step, every occurrence of a letter is marked and consecutive marked letters form a block.
A word is $k$-local if it can be marked
with at most $k$  blocks of marked letters
in every step.
For example, the word $\mathtt{reappear}$ is $2$-local witnessed, for instance, by the marking sequence $(\mathtt{e},\mathtt{a},\mathtt{r},\mathtt{p})$: marking the letters step by step leads to the marked words $\mathtt{r\overline{e}app\overline{e}ar}$,
$\mathtt{r\overline{ea}pp\overline{ea}r}$, $\mathtt{\overline{rea}pp\overline{ear}}$, $\mathtt{\overline{reappear}}$ (where the marked letters are visualised by overlines).
Some combinatorial properties of the locality have been investigated in \cite{Blocksequences}, e.g., relating the different
roles occurrences of letters in a word may play during the marking process.
In \cite{Connections}, the connection between the locality, cutwidth, and pathwidth of graphs was established by reductions showing that determining the locality of a word is NP-complete. Very recently, in \cite{klocalgraphs}, locality was generalised to coloured graphs.

\begin{figure}
    \centering
    \begin{tikzpicture}[scale=0.55]
        \node[shape=circle,draw=black] (a) at (-2,0) {$\ta$};
        \node[shape=circle,draw=black] (n) at (0,0) {$\tn$};
        \node[shape=circle,draw=black] (l) at (2,0) {$\tl$};
        \node[shape=circle,draw=black] (o) at (0,-2) {$\to$};
        \node[shape=circle,draw=black] (g) at (2,-2) {$\tg$};

        \path [-] (a) edge node[left] {} (n);
        \path [-] (l) edge node[left] {} (n);
        \path [-] (o) edge node[left] {} (n);
        \path [-] (g) edge node[left] {} (n);
        \path [-] (l) edge node[left] {} (o);
        \path [-] (l) edge node[left] {} (g);
        \path [-] (o) edge node[left] {} (g);

    \end{tikzpicture}
    \caption{Graph represented by the word \texttt{analog}.}
    \label{wordrepexample}
\end{figure}
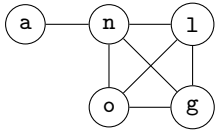

As started in \cite{wordsperspective}, we are going to connect this measure on patterns and words with the field of word-representable graphs.
A word {\em represents} a graph if the pairs of alternating letters exactly correspond
to the edges in the graph.
Two distinct letters are called {\em alternating} in a word if removing all other letters
does not yield two consecutive occurrences of the same letter.
Thus, the word's alphabet (the distinct letters occurring in it) determines the set of nodes and the projections onto two letters give the edge set.
For instance, the word $\mathtt{analog}$ represents the graph in \Cref{wordrepexample}.

The theory of word-representable graphs was introduced by Kitaev and
Pyatkin in \cite{onrep} as a tool for semigroup theory,
and an introduction to this theory can be found in \cite{wordsandgraphs}.
It has applications in periodic scheduling \cite{DBLP:conf/dlt/HalldorssonKP10,wordsandgraphs}, topology \cite{DBLP:journals/corr/abs-2111-10038}, and the power domination problem from physics \cite{chandrasekaran2019k},
but is also important from a theoretical perspective
since word-representable graphs generalise several important classes of graphs such as
circle graphs, 3-colourable graphs, and comparability graphs,
see \cite{wordsandgraphs}.
A fundamental result was presented in \cite{Semi}.
The authors characterised word-representable graphs
using a special type of graph orientation. 
The study of word-representable graphs focuses on 
word-representability of certain graphs classes,
operations that preserve word-representability,
the concept of $k$-representability,
and various generalisations of word-representability.
For instance, in \cite{Generalized}, word-representability
is generalised from the perspective of language theory.
The important notion of $k$-representability takes into account how many occurrences of a letter are necessary in order to represent the graph, i.e., a graph is $k$-representable for $k\in \N$
if it is represented by a word with at most $k$ copies of each letter. Thus, each complete graph on $n\in\N$ nodes $V=\{v_1,\ldots,v_n\}$ is $1$-representable witnessed by the word $v_1\cdots v_n$.
We denote the class of $k$-representable graphs by $\R^k$.
In \cite{wordsperspective}, the authors showed that $k$-representability
is related to the locality of words.
They proved that $\R^k$ is a subset of $\L^{2k}$ for every $k\in \N$,
where $\L^k$ denotes the class of graphs represented by $k$-local words.
There also $\L^1$ was characterised.

\textbf{Our Contribution.}
We further investigate the relation between locality and word-representability.
In particular, we prove that $\L^k$ is a strict subset of $\R^{k+1}$ for every $k\in \N$.
We show that $\L^1$ is the class of threshold graphs,
which means the classes $\L^k$ are a generalisation of them motivating further research.
Moreover,
we investigate the number of graphs in $\R^k$ and $\L^k$
and prove that $\L^k$ belongs to the factorial layer for $k>1$.
Finally, we show that the clique-width of graphs in $\L^k$ is bounded by $2^k+1$.

\textbf{Structure of the Work.}
After introducing the necessary definitions in \cref{prelims},
we prove that $\L^k$ is a subset of $\R^{k+1}$ in \cref{k+1}.
In \cref{speed}, we investigate the number of graphs
in $\R^k$ and $\L^k$.
We present our results regarding the clique-width in \cref{cwd}.
At the end, there is a conclusion in \cref{conclusion}.


%% file: prelims.tex
Let $\N$ denote the natural numbers starting with $1$ and set $\N_0=\N\cup\{0\}$
as well as $[n]=\{1,\ldots,n\}$ and $[n]_0=[n] \cup \{0\}$ for $n\in\N$.
For a set $M$, define $\binom{M}{2}=\{N\subseteq M \mid |N|=2\}$ as the set of all of $M$'s subsets with two elements.
Before we introduce our main definitions about word-representable graphs,
we present some basic definitions from combinatorics on words and graph theory.

An \emph{alphabet} $\Sigma$ is a finite set whose elements are called \emph{letters}.
A \emph{word} $w=\ta_1\cdots \ta_n$ is a finite sequence of letters $\ta_i\in\Sigma$ for $i\in[n]$,
$n\in\N_0$ from a given alphabet $\Sigma$ of length $n$, which is denoted by $|w|$.
The \nth{$i$} letter of $w$, i.e., $\ta_i$, is also denoted by $w[i]$ for all $i\in[|w|]$ if $|w|>0$.
The set of all finite words (strings, sequences) over the alphabet $\Sigma$ is denoted by $\Sigma^{\ast}$.
The empty word $\varepsilon$ is the word of length $0$.
Let $\Sigma^n$ denote all words in $\Sigma^{\ast}$ exactly of length $n\in\N_0$.
For all $\ta\in\Sigma$, define $|w|_{\ta}= |\{i\in[|w|] \mid w[i]=\ta\}|$ as the number of occurrences of $\ta$ in the word $w$.
Set $\alph(w) = \{\ta \in \Sigma \mid \exists i \in [|w|]: w[i] = \ta \}$ as $w$'s alphabet. 
A word $w$ is {\em $k$-uniform} for $k \in \N$ if $|w|_{\ta}=k$ for every $\ta \in \alph(w)$.
For $u,w\in\Sigma^{\ast}$, $u$ is called a \emph{factor}
of $w$, if $w = xuy$ for some words $x,y\in\Sigma^{\ast}$.
A mapping $f:\Sigma^{\ast}\rightarrow\Sigma^{\ast}$
is called a {\em morphism} if $f(uv)=f(u)f(v)$ holds for all $u,v\in\Sigma^{\ast}$.
Thus, a morphic mapping is already completely defined, if the images for all letters in $\Sigma$ are given.
For $S \subseteq \Sigma$, we define the {\em projection $\pi_S$ onto $S$} as the morphism
given by $\pi_S(x) = x$ if $x \in S$ and $\pi_S(x) = \varepsilon$ otherwise for all $x\in \Sigma$. 
For further definitions on words, see \cite{lothaire}.
This leads us to the main definition of the field of combinatorics on words for the word-representability.

\begin{definition}{}
Two distinct letters $\ta,\tb\in\Sigma$ {\em alternate} in a word $w\in\Sigma^{\ast}$
if $\pi_{\{\ta,\tb\}}(w)$ is $(\ta\tb)^n$, $(\ta\tb)^n\ta$, $(\tb\ta)^n$, or $(\tb\ta)^n\tb$ for some $n \in \N$.
\end{definition}

A \emph{finite, undirected graph} $G$ is a pair $G=(V,E)$ with finite sets $V$ of {\em nodes} and {\em edges} $E \subseteq \binom{V}{2}$.
By convention, set $n=|V|$ and $m=|E|$, and since we are only investigating finite, undirected graphs, we call them simply {\em graphs}. 
Moreover, given a graph $G$, let $V(G)$ and $E(G)$ denote its sets of nodes and edges resp.
A \emph{node labelled graph} is a graph $(V,E,f)$ with a function $f:V \rightarrow M$ for a set of {\em labels} $M$.
For $v,v'\in V$, a {\em path} from $v$ to $v'$ is a sequence $(v_1,\ldots,v_k)$
for some $k\in\N$ with $v=v_1$, $v'=v_k$, and $\{v_i,v_{i+1}\}\in E$ for all $i\in[k-1]$.
The {\em connected component of $v\in V$} is defined as $C(v)=\{u\in V \mid \exists \mbox{ path from } u \mbox{ to } v \}$.
A graph $G$ is called {\em connected} if $V=C(v)$ for all $v\in V$.
For $n \in \N$, the \emph{complete graph} on $|V| = n$ nodes is defined by $K_n = (V, \binom{V}{2})$
and the empty graph as $E_n=(V,\emptyset)$.
A graph $H = (U,F)$  is called a \emph{subgraph} of $G$ ($H \leq G$) if we have $U\subseteq V$ and $F\subseteq E\cap\binom{U}{2}$.
Moreover, we say that $H\leq G$ is \emph{induced} by $U$ if $F = E\cap \binom{U}{2}$; this is notated as $G|_{U}$.
A complete subgraph of $G$ is called a \emph{clique} and an empty subgraph {\em independent set}.
A {\em component} is a maximal connected induced subgraph.
For $n\in \N$, the {\em crown graph} $H_{n,n}$ is the graph $([2n],E)$ with $E=\{\{x,y+n\} \mid \text{$x,y\in [n]$ and  $x\neq y$}\}$.

With these definitions at hand, we define the main notion of interest, the word-representable graphs.

\begin{definition}
A graph $G = (V, E)$ is {\em represented} by a word $w \in \Sigma^*$ if $\alph(w)=V$
and for all $\ta, \tb \in V$, $\ta$ and $\tb$ alternate in $w$ iff $\{\ta, \tb\} \in E$.
A graph is {\em word-representable} if there is a word representing it.
A graph is {\em $k$-representable} for some $k \in \N$ if there exists a $k$-uniform word representing it.
The graph represented by a word $w$ is denoted by $\G(w)$.
\end{definition}

For example, the graph $G = (\{\mathtt{a,b,l,n,o}\}, \{\{\mathtt{b,a}\}, \{\mathtt{a,n}\},\{\mathtt{b,n}\}\})$ is represented by the word $\mathtt{balloon}$.

\begin{remark}\label{krep}
As shown in \cite{onrep}, a graph is $k$-representable iff there is a word $w$ representing it with $|w|_{\ta}\leq k$ for all $\ta \in \alph(w)$.
\end{remark}

We finish this section with our second notion of interest, the $k$-local words.
In the following, we give the main definitions on $k$-locality, following~\cite{LocalPatterns}. 

\begin{definition}
Let $\Sigma$ be an alphabet with $|\Sigma|=n\in\N$.
Set $\overline{\Sigma}=\{\overline{\ta} \mid \ta\in\Sigma\}$ as the alphabet containing the {\em marked letters} and
$\overline{w}=\overline{w[1]}\cdots \overline{w[|w|]}$ for $w\in \Sigma^*$.
We assume w.l.o.g. $\alph(w)=\Sigma$.
A {\em marking sequence} $\sigma$ is an enumeration of $\Sigma$, i.e.,
$\sigma=(\ta_1,\ldots,\ta_n)$ such that $\ta_i\neq\ta_j$ for all $i,j\in[n]$ with $i\neq j$. Let $i\leq n$.
We call $w_i$ the {\em marked version} of $w$ at stage $i$ w.r.t. $\sigma$ if $w_i[j] = \overline{w[j]}$ for all $j\in[|w|]$
for which we have $w[j] = \ta_\ell$ and $\ell \leq i$ with $\ell\in[n]$, and otherwise $w_i[j] = w[j]$.
A {\em marked block} (or just {\em block}) is a maximal consecutive factor of marked letters.
\end{definition}
\begin{definition}
A word $w\in\Sigma^{\ast}$ is called {\em $k$-local} for some $k\in\N$ if there is a marking sequence $\sigma$
such that $w_i$ contains at most $k$ marked blocks at every stage $i$.
It is called {\em strictly $k$-local} if it is $k$-local but not $(k-1)$-local.
\end{definition}

For example, consider the word $\mathtt{pepper}$ and the marking sequence $\sigma = (\mathtt{r}, \mathtt{p},\mathtt{e})$.
We have $\t{peppe}\overline{\tr}$ at the first stage,
$\overline{\tp}\te\overline{\tp\tp}\te\overline{\tr}$ at the second stage, 
and $\overline{\t{pepper}}$ at the third stage.
At every stage, there are at most three blocks, which means that this word is $3$-local.
One can verify that it is strictly $2$-local witnessed by the marking sequence $\sigma' = (\mathtt{p}, \mathtt{e},\mathtt{r})$.

\begin{definition}
We denote the class of word-representable graphs by $\R$.
For every $k \in \N$, we denote the class of $k$-representable graphs by $\R^k$ and the class of graphs representable by $k$-local words by $\L^k$.
\end{definition}

We present some section specific definitions within the following sections.

%% file: k+1.tex
In \cite{wordsperspective}, it has been shown that every $k$-local word for $k\in \N$ represents a $2k$-representable graph.
We improve this result by showing that every $k$-local word also represents a $(k+1)$-representable graph.
To do so, we need the following two lemmata.

\input{./Statements/lemmaaltk+1}
\ifpaper
\else
\input{./Proofs/proof_lemmaaltk+1}
\fi

\input{./Statements/lemmastillkloc}
\ifpaper
\else
\input{./Proofs/proof_lemmastillkloc}
\fi

Using the previous two lemmata, for each $k$-local word $w$,
we can provide a construction such that we obtain a $k$-local and $(k+1)$-uniform word $w'$ that represents the same graph.

\input{./Statements/lemmalemk+1}
\ifpaper
\else
\input{./Proofs/proof_lemmalemk+1}
\fi

\begin{remark}
Notice that if $w$ is $k$-local for some $k\in\N$ there exists a $k'\leq k$ such that $w$ is strictly $k'$-local.
Since $w$ remains the same, we can alter \Cref{lemk+1} to the existence of
a $k'$-local word $w'$ with $|w'|_{\ta}\leq k'+1$ for every $\ta\in\alph(w')$ such that $\G(w)=\G(w')$. 
\end{remark}

\Cref{lemk+1}, together with \Cref{krep}, immediately implies the following result.

\begin{theorem}\label{repk+1}
Every $k$-local word for $k\in \N$ represents a $(k+1)$-representable graph.
\end{theorem}

We have proven $\L^k \subseteq \R^{k+1}$ for every $k\in \N$, but we do not know yet if the inclusion is strict.
In \cref{cwd}, we will show, among other results, that this is the case.



%% file: Statements/lemmaaltk+1.tex
\begin{lemma}{}\label{altk+1}
Let $w \in \Sigma^*$ be a $k$-local word for $k \in \N$.
Every letter $\ta \in \alph(w)$ with $|w|_{\ta}>k+1$
does not alternate with any other letter.
\end{lemma}

%% file: Proofs/proof_lemmaaltk+1.tex
\begin{proof}
Let $\ta \in \alph(w)$ with $|w|_{\ta}>k+1$.
Let $\tb \in \alph(w)\m \{\ta\}$.
Let $s$ be a marking sequence of $w$ witnessing its $k$-locality.
For every $p \in [|\alph(w)|]$, $w_p$ w.r.t. $s$ is of the form
$$u_1\overline{m_1}u_2\overline{m_2}\cdots u_k\overline{m_k}u_{k+1}$$
for some $u_1,m_1,\ldots,u_k,m_k,m_{k+1} \in \Sigma^*$.

First, assume $\tb$ occurs before $\ta$ in $s$ at position $p$.
This implies $|m_i|_{\ta}=0$ for every $i\in [k]$
and $|u_j|_{\tb}=0$ for every $j \in [k+1]$.
Since $|w|_{\ta}>k+1$, there is $j \in [k+1]$ with $|u_j|_{\ta}>1$.

Now, assume $\ta$ occurs before $\tb$ in $s$ at position $p$.
This implies $|m_i|_{\tb}=0$ for every $i\in [k]$
and $|u_j|_{\ta}=0$ for every $j \in [k+1]$.
Since $|w|_{\ta}>k+1$, there is $i \in [k]$ with $|m_i|_{\ta}>1$.

In both cases, there is a factor of $w$ with two occurrences of $\ta$ and no occurrence of $\tb$.
This means $\ta$ does not alternate with $\tb$ in $w$.\qed
\end{proof}

%% file: Statements/lemmastillkloc.tex
\begin{lemma}{}\label{stillkloc}
    Let $w \in \Sigma^*$ be $k$-local and $S\subseteq \alph(w)$.
    Then $\pi_S(w)$ is also $k$-local.
\end{lemma}

%% file: Proofs/proof_lemmastillkloc.tex
\begin{proof}
    Let $s$ be a marking sequence witnessing the $k$-locality.
    Let $i \in [n]$.
    At every stage, the marked version of $w$ w.r.t. $s$ is of the form
    $$u_1\overline{m_1}u_2\overline{m_2}\cdots u_k\overline{m_k}u_{k+1}$$
    for some $u_1,m_1,\ldots,u_k,m_k,m_{k+1} \in \Sigma^*$,
    and the marked version of $\pi_S(w)$ w.r.t. $s$ is
    $$\pi_S(u_1)\overline{\pi_S(m_1)}\pi_S(u_2)\cdots \pi_S(u_{k})\overline{\pi_S(m_{k})}\pi_S(u_{k+1}).$$
    Since there are at most $k$ blocks, $\pi_S(w)$ is $k$-local.\qed
\end{proof}

%% file: Statements/lemmalemk+1.tex
\begin{lemma}\label{lemk+1}
Let $k \in \N$. Let $w$ be a $k$-local word.
There is a $k$-local word $w'$
with $|w'|_{\ta}\leq k+1$ for every $\ta \in \alph(w')$ such that $\G(w) = \G(w')$.
\end{lemma}

%% file: Proofs/proof_lemmalemk+1.tex
\begin{proof}
Let $S=\{\ta_1,\ldots \ta_\ell\}$ be the set of letters that occur more than $k+1$ times in $w$.
Define
$$w'=\pi_{\alph(w)\m S}(w)\ta_1\ta_1\cdots \ta_\ell\ta_\ell.$$
We have $|w'|_{\ta}\leq k+1$ for every $\ta \in \alph(w')$.
By \Cref{altk+1}, $w$ and $w'$ represent the same graph.
By \cref{stillkloc}, $\pi_{\alph(w)\m S}(w)$ is $k$-local.
Since we can mark $\ta_1,\ldots,\ta_\ell$ at the end in this order,
$w'$ is also $k$-local.\qed
\end{proof}

%% file: speed.tex
In this section we first show that $\R,\R^k$, and $\L^k$ are hereditary. This leads to the first main result of this section, namely that $\L^1$ is exactly the class of threshold graphs.
The second main result stems from the fact that hereditary graph classes have an asymptotic growth in their size w.r.t. their number of nodes.
For this purpose, we use the notion of {\em speed} and give an upper and a lower bound for $\R^k$ and $\L^k$.

\begin{definition}
A class of graphs is called {\em hereditary} if it is closed under induced subgraphs.
By $\Free(M)$, we denote the class of graphs that do not contain a graph from the set of graphs $M$ as an induced subgraph.
The graphs from $M$ are called the {\em forbidden induced subgraphs} of $\Free(M)$.
A class of graphs $X$ is called {\em finitely defined} if there is a finite set $M$ such that $X=\Free(M)$.
\end{definition}

\begin{remark}
Note that for a finitely defined class, it can be checked in polynomial time in the size of a graph if it is contained in the class.
A class of graphs $X$ is hereditary iff $X = \Free(M)$ for some set of graphs $M$ (cf. \cite{wordsandgraphs}).
\end{remark}

Therefore, every hereditary class of graphs can be characterised in terms of forbidden induced subgraphs,
and is a crucial task to find a characterisation with a minimal set of forbidden induced subgraphs.
However, for many important graph classes such a characterisation has not yet been found.
First, we prove that the classes $\R$, $\R^k$, $\L^k$ are hereditary.

\input{./Statements/proprlhereditary}

\input{./Proofs/proof_proprlhereditary}

We provide a forbidden induced subgraph characterisation for $\L^1$ using the following theorem. Notice that we hereby connect the notion of $1$-local representable graphs with the one of threshold graphs \cite{isgci}.

\begin{definition}
	A graph $G$ is a \emph{threshold graph} if it can be constructed from the empty graph by repeatedly adding either an isolated node
	or a node that is connected to all other nodes.
\end{definition}

\begin{theorem}
$\L^1$ is exactly the class of threshold graphs.
\end{theorem}
\begin{proof}
In \cite{wordsperspective}, the authors introduced the class of $1$-local representable graphs 
and showed that it contains exactly the graphs represented by $1$-local words, i.e., $\L^1$.
The definition of threshold graphs and $1$-local representable graphs is the same.
\qed
\end{proof}

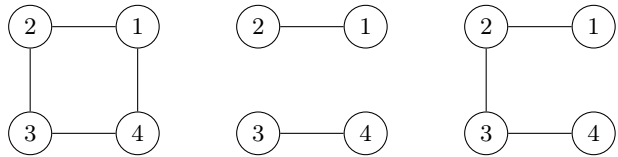
\begin{figure}
    \begin{tikzpicture}[scale = 1]
    \begin{scope}[shift={(0,0)}]
          \node[shape=circle,draw=black] (1) at (45:1) {1};
          \node[shape=circle,draw=black] (2) at (135:1) {2};
          \node[shape=circle,draw=black] (3) at (225:1) {3};
          \node[shape=circle,draw=black] (4) at (315:1) {4};
  
          \path [-](1) edge node {} (2);
          \path [-](2) edge node {} (3);
          \path [-](3) edge node {} (4);
          \path [-](4) edge node {} (1);
    \end{scope}
    \begin{scope}[shift={(3,0)}]
        \node[shape=circle,draw=black] (1) at (45:1) {1};
        \node[shape=circle,draw=black] (2) at (135:1) {2};
        \node[shape=circle,draw=black] (3) at (225:1) {3};
        \node[shape=circle,draw=black] (4) at (315:1) {4};

        \path [-](1) edge node {} (2);
        \path [-](3) edge node {} (4);
    \end{scope}
    \begin{scope}[shift={(6,0)}]
        \node[shape=circle,draw=black] (1) at (45:1) {1};
        \node[shape=circle,draw=black] (2) at (135:1) {2};
        \node[shape=circle,draw=black] (3) at (225:1) {3};
        \node[shape=circle,draw=black] (4) at (315:1) {4};

        \path [-](1) edge node {} (2);
        \path [-](2) edge node {} (3);
        \path [-](3) edge node {} (4);
    \end{scope}
    \end{tikzpicture}
    \centering
    \caption{Forbidden induced subgraphs of $\L^1$.}
    \label{forbidden}
\end{figure}

It is well known that a graph is a threshold graph iff it does not contain one of the forbidden induced subgraphs
shown in Figure~\ref{forbidden}, cf. \cite{threshold}.
Note that the class of threshold graphs is equivalent to several other classes,
which are listed in \cite{isgci}.

This means that the classes $\mathcal L^k$ for $k\in \N$ are a generalisation of the important threshold graphs,
which further motivates to investigate them.
However, we have not found a forbidden induced subgraph characterisation of $\L^k$ for $k>1$, 
and it is unknown if these classes are finitely defined.

Another important aspect of hereditary graph classes is the asymptotic growth of their size w.r.t. the number of nodes.
First, we define a way to count the size in our case (cf. \cite{speedrep}).

\begin{definition}
For a class of graphs $X$ and $n \in \N$ we define
$X_n=\{G \in X \mid V(G)=[n]\},$
i.e. the graphs in $X$ with nodes from $1$ to $n$.
We call $n \mapsto |X_n|$ the \emph{speed} of $X$.
\end{definition}

For $X=\R$, the speed was analysed in \cite{speedrep}.
The authors proved $|\R_n|=2^{\frac{n^2}{3}+o(n^2)}$ and concluded that there are
$2^{\frac{n^2}{3}+o(n^2)}$ non-isomorphic word-representable
graphs with $n$ nodes.
We will investigate the speed of $\L^k$ and $\R^k$ with a similar approach based on the index (also called \emph{colouring number}) of a hereditary class.
See \cite{HereditaryAndMonotone} for an introduction.

\begin{definition}
Let $\E_{i,j}$ be the class of graphs whose
nodes can be partitioned into $i$ independent sets and $j$ cliques as induced subgraphs.
The \emph{index} $c(X)$ of a hereditary class of graphs $X$ is
the largest $k \in \N$ such that
$\E_{i,k-i} \subseteq X$ holds for some $i \in \N$.
\end{definition}

We want to determine the index of $\R^k$ and $\L^k$ for every $k \in \N$.
First, we observe that there are graphs in $\E_{0,2}$ and $\E_{1,1}$ that are not word-representable.
These graphs are depicted in Figure~\ref{G1} and Figure~\ref{G2}.
\input{./Statements/lemE02E11}
\ifpaper
\else
\input{./Proofs/proof_lemE02E11}
\fi

\begin{figure}
\centering
\begin{minipage}{.4\textwidth}
	\begin{tikzpicture}[scale = 0.8]
		\node[shape=circle,draw=black,fill=orange] (3) at (90:2) {3};
		\node[shape=circle,draw=black,fill=orange] (4) at (162:2) {4};
		\node[shape=circle,draw=black,fill=cyan] (5) at (234:2) {5};
		\node[shape=circle,draw=black,fill=cyan] (6) at (306:2) {6};
		\node[shape=circle,draw=black,fill=orange] (2) at (18:2) {2};
		\node[shape=circle,draw=black,fill=orange] (1) at (0:0) {1};
        \node[shape=circle,draw=black,fill=cyan] (7) at (270:1) {7};
	
		\path [-](3) edge node {} (4);
		\path [-](4) edge node {} (5);
		\path [-](5) edge node {} (6);
		\path [-](6) edge node {} (2);
		\path [-](2) edge node {} (3);
		\path [-](1) edge node {} (3);
		\path [-](1) edge node {} (4);
		\path [-](1) edge node {} (5);
		\path [-](1) edge node {} (6);
		\path [-](1) edge node {} (2);
        \path [-](7) edge node {} (2);
        \path [-](7) edge node {} (4);
        \path [-](7) edge node {} (5);
  		\path [-](7) edge node {} (6);
        \path [-](2) edge node {} (4);
	\end{tikzpicture}
	\centering
    \caption{A graph in $\E_{0,2}$ that is not in $\R$ (the colouring indicates a partition into two cliques).}
    \label{G1}
\end{minipage}
\hfill
\begin{minipage}{.4\textwidth}
	\begin{tikzpicture}[scale = 0.8]
		\begin{scope}[shift={(0,0)}]
			\node[shape=circle,draw=black,fill=orange] (1) at (0:0) {1};
			\node[shape=circle,draw=black,fill=orange] (2) at (30:1) {2};
			\node[shape=circle,draw=black,fill=orange] (3) at (150:1) {3};
			\node[shape=circle,draw=black,fill=orange] (4) at (270:1) {4};
			\node[shape=circle,draw=black,fill=cyan] (5) at (90:2) {5};
			\node[shape=circle,draw=black,fill=cyan] (6) at (210:2) {6};
			\node[shape=circle,draw=black,fill=cyan] (7) at (330:2) {7};
			
			\path [-](1) edge node {} (2);
			\path [-](1) edge node {} (3);
			\path [-](1) edge node {} (4);
			\path [-](2) edge node {} (3);
			\path [-](3) edge node {} (4);
			\path [-](4) edge node {} (2);
			\path [-](5) edge node {} (2);
			\path [-](5) edge node {} (3);
			\path [-](6) edge node {} (3);
			\path [-](6) edge node {} (4);
			\path [-](7) edge node {} (2);
			\path [-](7) edge node {} (4);
		\end{scope}
	\end{tikzpicture}    
	\centering
	\caption{A graph in $\E_{1,1}$ that is not in $\R$ (the colouring indicates a partition into a clique and an independent set).}
	\label{G2}
\end{minipage}
%
\end{figure}

With this we can show that the indices of $\R^{k+1}$ and $\L^k$ are both $1$ for every $k\in\N$.

\begin{theorem}
For every $k \in \N$, $c(\R^{k+1})=c(\L^k)=1$ holds.
\end{theorem}
\begin{proof}
Let $k \in \N$.
Note that $\E_{1,0}$ contains only the empty graphs, all of which are in $\R^{k+1}$ and $\L^k$.
Due to \cref{E02E11}, we only need to show that $\E_{2,0}$ is not contained in $\R^{k+1}$ and $\L^k$.
We know from \cite{crown} that the crown graph $H_{2k+3,2k+3}$ has representation number $\ceil{k+\frac32}>k+1$.
By \cref{repk+1}, this implies that it is neither contained in $\R^{k+1}$ nor in $\L^k$.
Since this graph is $2$-colourable and thus in $\E_{2,0}$, it follows that $\E_{2,0}$ is not contained in $\R^{k+1}$ and $\L^k$.
Hence, we have $c(\R^{k+1})=c(\L^k)=1$.\qed
\end{proof}

The speed is closely related to the entropy, which can be considered as a coefficient of compressibility (cf. \cite{wordsandgraphs}).
\begin{definition}
    For a class of graphs $X$, its \emph{entropy} is defined by
    $$e(X)=\lim_{n\rightarrow\infty}\frac{\log_2 |X_n|}{\binom{n}{2}}.$$
\end{definition}

With the Alekseev–Bollobás–Thomason Theorem (cf. \cite{Alekseev}) we immediately get the following result. 
\begin{corollary}
For $X \in \{\R^{k+1},\L^k\}$ with $k \in \N$,
the entropy is $e(X)=0$.
\end{corollary}

Unfortunately, this only gives us an upper bound for the speed of $\R^k$ and $\L^k$.
We continue with a lower bound based on the Bell number.

\begin{definition}
For every $n \in \N$, the \nth{$n$} \emph{Bell number} $B_n$ is defined as the number of partitions of a set with $n$ elements.
\end{definition}

The hereditary properties are partitioned into properties below the Bell number with $|X_n|<B_n$ and properties above the Bell number with $|X_n|\geq B_n$ for sufficiently large $n\in \N$.
The properties below the Bell number are of particular interest since they are finitely defined (\cite{warwick67023})
and all contained graphs have bounded clique-width (\cite{cliquewidth}).

\ifpaper
\input{./Statements/lemcliques}
\else
\input{./Statements/lemcliques}
\input{./Proofs/proof_lemcliques}
\fi

\begin{theorem}{}\label{AboveBell}
Let $X \in \{\R^{k},\L^{k}\}$ for $k>1$.
We have $B_n\leq |X_n|$.
\end{theorem}
\begin{proof}
Let $n\in\N$, and let $P$ be a partition of $[n]$.
Set $G=([n],E)$ with
$$E=\bigcup_{C\in P}\{\{x,y\}\mid x,y \in C\}.$$
By \Cref{cliques}, $G \in X_n$.
Since we get a unique graph for every partition, it follows $B_n\leq |X_n|$.\qed
\end{proof}

Thus, $\L^k$ and $\R^k$ are above the Bell number for $k>1$. 
This means that we cannot draw any conclusion regarding forbidden induced subgraphs or the clique-width of these classes.
Nevertheless, we are going to inspect the clique-width of $\L^k$ from a constructive perspective.

%% file: Statements/proprlhereditary.tex
\begin{theorem}\label{rlhereditary}
The classes $\R$, $\R^k$, and $\L^k$ are hereditary for all $k\in\N$.
\end{theorem}

%% file: Proofs/proof_proprlhereditary.tex
\begin{proof}
If we have a word $w$ and remove all occurrences of a letter $\ta$,
the resulting word represents the graph $\G(w)|_{\alph(w)\m\{\ta\}}$.
If $w$ is $k$-uniform or $k$-local for some $k\in \N$,
the resulting word is still $k$-uniform or $k$-local, respectively.
Therefore, $\R$, $\R^k$, and $\L^k$ are hereditary for every $k \in \N$.\qed
\end{proof}

%% file: Statements/lemE02E11.tex
\begin{lemma}{}\label{E02E11}
We have $\E_{0,2} \not\subseteq \R$ and $\E_{1,1} \not\subseteq \R$.
\end{lemma}

%% file: Proofs/proof_lemE02E11.tex
\begin{proof}
The graph depicted in Figure~\ref{G1} can be partitioned into two cliques as indicated by the colouring.
Therefore, it belongs to $\E_{0,2}$.
Since we know from \cite{humanverifiable} that this graph is not word-representable,
$\E_{0,2}$ is not a subset of $\R$.

The graph depicted in Figure~\ref{G2} can be partitioned into a clique and an independent set as indicated by the colouring.
Therefore, it belongs to $\E_{1,1}$.
Since we know from \cite{humanverifiable} that this graph is not word-representable,
$\E_{1,1}$ is not a subset of $\R$.\qed
\end{proof}

%% file: Statements/lemcliques.tex
\begin{lemma}\label{cliques}
Let $G$ be a graph with every component being a clique.
Then, $G \in X$ for $X \in \{\R^{k},\L^{k}\}$ for $k>1$.
\end{lemma}

%% file: Proofs/proof_lemcliques.tex
\begin{proof}
Let $C_1,\ldots,C_n$ be the components of $G$.
Let $w_i$ be a word containing every node from $C_i$ exactly once for every $i \in [n]$.
Then, $w= w_1w_1w_2w_2 \cdots w_nw_n$ represents $G$.
Note that $w$ is $2$-uniform and $2$-local, witnessed by the marking sequence
$$(w_1[1],\ldots,w_1[|w_1|],\ldots,w_n[1],\ldots,w_n[|w_n|]).$$
This implies $G\in \R^k$ and $G\in \L^k$.\qed
\end{proof}

%% file: cwd.tex
In this last section, we prove that $\L^k$ is a proper subset of $\R^{k+1}$ for all $k\in\N$.
Therefore, we first introduce the well-known graph parameter {\em clique-width} and show
that $\L^k$ has bounded clique-width based on $k$.
Afterwards, we show that $\L^k$ belongs to the factorial layer, a subclass of the hereditary graph classes depending on their speed.

We start with a folklore definition on constructing labelled graphs.
Note that every labelled graph can be constructed using only these operations.

\begin{definition}
We define the following operations on labelled graphs:\\
- $i(v)$, create a node $v$ with label $i$,\\
- $G \oplus H$, disjoint union of $G$ and $H$,\\
- $\eta_{i,j}$, connect every node with label $i$ to every node with label $j$ for $i \neq j$, \\
- $\rho_{i \rightarrow j}$, rename label $i$ to $j$.

\medskip

A \emph{$k$-expression} with $k\in\N$ is an expression consisting of these operations with not more than $k$ distinct labels.
For a graph $G$, its {\em clique-width} $\cwd(G)$ is the minimum number such that there is $k$-expression for a labelled version of it.
We define an \emph{$M$-expression} for a set $M$ as an expression with the above operations using only labels from $M$.
For every $k\in \N$, we define $\Sigma_k=\{0,1\}^k\cup \{2\}$.
For a marked word $m\in(\Sigma\cup\overline{\Sigma})^*$
and $k\in \N$, define the function $\lab_m^k: \Sigma \rightarrow \Sigma_k$ that maps a letter $\ta$ to a tuple
such that the \nth{$i$} component is $1$ if there is one occurrence of $\ta$ in the \nth{$i$} block of $m$
and it is $0$ if there are less than $i$ blocks or there is no occurrence of $\ta$ in the \nth{$i$} block.
If there is a block with two or more occurrences, $\ta$ is mapped to $2$ instead.
\end{definition}

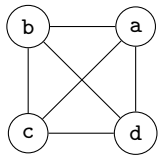
\begin{figure}
    \begin{tikzpicture}[scale = 1]
    \begin{scope}[shift={(0,0)}]
          \node[shape=circle,draw=black] (1) at (45:1) {$\ta$};
          \node[shape=circle,draw=black] (2) at (135:1) {$\tb$};
          \node[shape=circle,draw=black] (3) at (225:1) {$\tc$};
          \node[shape=circle,draw=black] (4) at (315:1) {$\td$};
  
          \path [-](1) edge node {} (2);
          \path [-](2) edge node {} (3);
          \path [-](3) edge node {} (4);
          \path [-](4) edge node {} (1);
          \path [-](1) edge node {} (3);
          \path [-](4) edge node {} (2);
    \end{scope}
    \end{tikzpicture}
    \centering
    \caption{Graph with clique-width 2.}
    \label{cwdex}
\end{figure}

    The graph in Figure~\ref{cwdex} can be constructed using only the labels $1$ and $2$ with the expression
    $$\eta_{1,2}(\rho_{2 \rightarrow 1}(\eta_{1,2}(1(\ta) \oplus 2(\tb)))  \oplus  \rho_{1 \rightarrow 2}(\eta_{1,2}(1(\tc) \oplus 2(\td)))).$$
    Note that a graph with clique-width $1$ cannot have any edges because we can only add edges between distinct labels.
    Therefore, this graph has a clique-width of $2$.

Using these definitions, we can show that $\L^k$ has a bounded clique-width based on $k$.
For the proof, we need the following lemma.

\begin{lemma}{}\label{expalt}
Let $w$ be a $k$-local word witnessed by a marking sequence $s$.
Let $i\in [|s|]$ and $\ta = s[i]$.
For every $j \in [|s|]_0$, let $u_j=\pi_{\alph(s[1]\cdots s[j])}(w)$ and
let $w_j$ be the marked version of $w$ at stage $j$ w.r.t. $s$.
Let $E$ be a $\Sigma_k$-expression for $\G(u_{i-1})$ such that every node $v$ has label
$\lab_{w_{i-1}}^k(v)$.
There is a $\Sigma_k$-expression for $\G(u_i)$ such that every node $v$ has label $\lab_{w_{i}}^k(v)$.
\end{lemma}
\begin{proof}
Let $N$ be the tuple with $0$ in all $k$ components.
Let $l_1,\ldots,l_n$ for $n\in \N_0$ be the labels of the letters alternating with $\ta$ in $u_i$.
Note that for each letter alternating with $\ta$, every letter with the same label also alternates with $\ta$ since
they have the same number of occurrences in each block.
Thus,
$$\eta_{l_1,N}(\ldots\eta_{l_n,N}(N(\ta) \oplus E)\ldots)$$
is a $\Sigma_k$-expression for $\G(u_{i})$.
Because no label is changed, every node $v$ except $\ta$ has label $\lab_{w_{i-1}}^k(v)$.

Now, we need to adjust the labels because blocks can merge and new blocks can arise.
If the \nth{$j$} and the \nth{$(j+1)$} block merge, we need to add the numbers of occurrences in both blocks and append a $0$ at the end.
This means every label $(x_1,\ldots,x_p)$ is renamed to $(x_1,\ldots,x_{j-1},x_j+x_{j+1},x_{j+2},\ldots,x_p,0)$ if $x_j+x_{j+1}<2$, else to $2$,
and $2$ is not renamed.
We define a function $g$ that maps the labels of the letters from $u_{i-1}$ to the corresponding labels after all merges.
This function induces a directed graph with self-loops
where every label has an edge to the new label after all merges.
Note that there is no cycle except for self-loops in this graph.
Hence, we find a topological ordering $l_1,\ldots,l_n$.
Let $w'$ be the word that we get from $w_{i-1}$ by marking every occurrence of $\ta$ that does not create a new block.
It follows that
$$E'=\rho_{l_1\rightarrow g(l_1)}(\ldots\rho_{l_n\rightarrow g(l_n)}(E)\ldots)$$
is a $\Sigma_k$-expression for $\G(u_{i})$ such that every node $v$ has label $\lab_{w'}^k(v)$.

If there arises a new block between the \nth{$j$} and the \nth{$(j+1)$} one, we need to insert a component with value $0$ into the tuples.
This means every label $(x_1,\ldots,x_p)$ is renamed to $(x_1,\ldots,x_{j},0,x_{j+1},\ldots,x_{p-1})$ and $2$ is not renamed.
We define a function $h$ that maps the labels of the letters from $u_{i-1}$ to the corresponding labels after all new blocks are there.
Again, this induces a directed graph with self-loops.
Note that there is no cycle except for self-loops in this graph.
Hence, we find a topological ordering $k_1,\ldots,k_n$.
Also, $\ta$ has to get the correct label.
Let $y = \lab_{w_{i}}^k(\ta)$.
It follows that
$$\rho_{N\rightarrow y}(\rho_{k_1\rightarrow h(k_1)}(\ldots\rho_{k_n\rightarrow h(k_n)}(E')\ldots))$$
is a $\Sigma_k$-expression for $\G(u_{i})$ such that every node $v$ has label $\lab_{w_{i}}^k(v)$.\qed
\end{proof}

Now, we can prove the first main theorem of this section.

\begin{theorem}{}\label{cwdLk}
Let $G \in \L^k$ for $k \in \N$.
We have $\cwd(G)\leq 2^k+1$.
\end{theorem}
\begin{proof}
We apply \cref{expalt}.
Starting with the empty expression, we get a $\Sigma_k$-expression for $G$ by induction.
It follows $\cwd(G)\leq |\Sigma_k|=2^k+1$.\qed
\end{proof}

    Consider the $2$-local word $\t{banana}$ and the marking sequence $(\tn,\ta,\tb)$.
    We start with
    $$E_1 = \rho_{N\rightarrow (1,1)}(N(\tn)).$$
    Next, we add $\ta$ alternating with all nodes that have the label $(1,1)$ and get
    $$E_2 = \rho_{N\rightarrow 2}(\rho_{2 \rightarrow 2}(\rho_{(1,1)\rightarrow 2}(\eta_{(1,1),N}(N(\ta) \oplus E_1)))).$$
    Here, $\tn$ needs to get the label $2$ since the first and the second block merge such that there are two copies of $\tn$ in the first block.
    Finally, we have
    $$E_3 = \rho_{N\rightarrow (1,0)}(\rho_{2 \rightarrow 2}(\rho_{2\rightarrow 2}(N(\tb) \oplus E_2))),$$
    which is a $\Sigma_2$-expression for $\G(\t{banana})$.

\begin{remark}
The theorem implies that $\L^k$ has bounded clique-width for every $k\in \N$.
This is interesting from an algorithmic point of view since many problems that cannot be solved in polynomial time in general
can be solved in polynomial time on classes with bounded clique-width.
These problems include the Hamiltonian path problem and the minimum maximal matching problem, see \cite{NPinPol}.
Courcelle et al. showed in \cite{logic} that all problems which can be formulated in a certain form of monadic second-order logic called MSO\textsubscript{1}
can even be solved in linear time
on graphs with clique-width bounded by some $k\in\N$ given a $k$-expression.
An example for such a problem is 3-colourability, see \cite{3colourability}. 
Note that for a class with bounded clique-width,
there is $k\in \N$ such that a $k$-expression for every graph in this class can be determined in polynomial time by \cite{approximate}.
In \cite{cwdNP}, it was shown that it is NP-complete to decide if the clique-width of a graph $G$ is at most $k$.
There are several other graph parameters related to clique-width.
By \cite{treewidth}, a graph with treewidth $k$ has a clique-width not greater than $3\cdot 2^{k-1}$.
By \cite{NLC}, a graph with NLC-width $k$ has a clique-width not smaller than $k$ and not greater than $2k$.  
By \cite{approximate}, a graph with rank-width $k$ has a clique-width not smaller than $k$ and not greater than $2^{k+1}-1$.  A thorough investigations of these graph parameters for $\L^k$ remains for future work.
\end{remark}

With the following theorem by Allen et al., we can apply our knowledge about the clique-width of $\L^k$ to learn more about its speed.

\begin{theorem}[{{\cite{cliquewidth}}}] 
If $X$ is a class of graphs with bounded clique-width, there is $C$ such that $|X_n|\leq n!C^n$.
\end{theorem}

\begin{corollary}
For every $k\in \N$, there is a constant $C$ such that $|\L_n^k|\leq n!C^n$.
\end{corollary}

Hereditary classes of graphs are divided into layers by their speed, first introduced in \cite{size}.
One of these layers is the factorial layer.
A class with speed $f(n)$ belongs to the factorial layer if there are constants $C_1,C_2>0$ such that
$$n^{C_1n}\leq f(n)\leq n^{C_2n}.$$ 
Since the factorial layer is split into classes above and below the Bell number (\cite{DecidingTheBellNumber}),
$\L^k$ belongs to the factorial layer for every $k > 1$ by \cref{AboveBell}.

To finish off this section, we can use our knowledge about the clique-width of $\L^k$ to show that $\L^k$ is a proper subset of $\R^{k+1}$. 

\begin{theorem}\label{ProperSubset}
For every $k\in \N$,
$\L^k \subsetneq \R^{k+1}$.
\end{theorem}
\begin{proof}
By \cite{Unbounded}, the class of circle graphs has unbounded clique-width.
The circle graphs are exactly the $2$-representable graphs,
as shown in \cite{onrep}.
Thus, $\R^{k+1}$ has unbounded clique-width for every $k\in \N$.
By \cref{cwdLk} and \cref{repk+1}, we get
$\L^k \subsetneq \R^{k+1}$.\qed
\end{proof}

%% file: conclusion.tex
In this work, we investigated the connection between
locality of words and word-representable graphs.
We proved that $\L^k$ is a subset of $\R^{k+1}$
for every $k\in \N$, and we later proved that this inclusion is strict.
We found out that $\L^1$ is the class of threshold graphs,
which gave us a characterisation in terms of forbidden induced subgraphs.
Moreover, we analysed the speed of $\R^k$ and $\L^k$.
In particular, we showed that both classes are above the Bell number for $k>1$
and $\L^k$ belongs to the factorial layer.
Finally, we showed that the clique-width of graphs in $\L^k$
is bounded by $2^k+1$, which connects the complexity of words with the 
complexity of the graphs represented by them.

It remains an open problem to characterise $\L^k$ for $k>1$.
In particular, it is unknown if $\L^k$ is finitely defined for $k>1$.
It is unclear if $\R^k$ also belongs to the factorial layer.
Also, we do not know if $2^{k+1}$ is optimal as an upper bound of the clique-width of $\L^k$.
Moreover, a thorough investigation of well-known graph-parameters for $\L^k$ remains future work.